\documentclass[14pt,reqno,a4pape]{amsart}
\usepackage{times,mathrsfs,color,comment}
\usepackage{amssymb,amsmath,bbm}
\usepackage[notref,notcite]{showkeys}

\numberwithin{equation}{section}

\newtheorem{thm}{Theorem}[section]
\newtheorem{cor}[thm]{Corollary}
\newtheorem{lem}[thm]{Lemma}
\newtheorem{ex}[thm]{Example}
\newtheorem{pro}[thm]{Proposition}

\newtheorem{rmk}[thm]{Remark}

\newcommand{\re}{{\rm Re\,}}

\newcommand{\D}{{\mathbb D}}

\begin{document}

\title[]
{Difference of composition operators over Bergman spaces with exponential weights}

\author[I. Park]{Inyoung Park}
\address{Institute for Basic Sciences, Korea University, Seoul 02841, KOREA}
\email{iypark26@gmail.com}

\thanks{The author was supported by the National Research Foundation of Korea(NRF) grant funded by the Korea government(MSIT) (No. 2021R1I1A1A01047051).}


\begin{abstract}
In this paper, we obtain a complete characterization for the compact difference of two composition operators acting on Bergman spaces with weight $\omega=e^{-\eta}$, $\Delta\eta>0$ in terms of the $\eta$-derived pseudodistance of two analytic self maps. In addition, we provide simple inducing maps which support our main result. We also study the topological path component of the space of all bounded composition operators on $A^2(\omega)$ endowed with the Hilbert-Schmidt norm topology.
\end{abstract}

\maketitle

\section{Introduction}
Let $S(\D)$ be the space of holomorphic self-maps of the unit disk $\D$. Given $\varphi$ in $S(\D)$, the composition operator $C_\varphi$ is defined by $C_\varphi f:=f\circ\varphi$ for all $f$ belonged to the holomorphic function spaces $H(\D)$. For the integrable radial function $\omega$, let $L^p(\omega dA)$ be the space of all measurable functions $f$ on $\D$ such that
\begin{align*}
\|f\|^p_{p}:=\int_\D |f(z)\omega(z)|^pdA(z)<\infty,\quad0<p<\infty,
\end{align*}
where $dA(z)$ is the normalized area measure on $\D$. We put $\|\cdot\|:=\|\cdot\|_2$ for simplicity. Especially, we denote $A^p(\omega):=L^p(\omega dA)\cap H(\D)$ and we use the notation $A^p_\alpha(\D)$ when $\omega(z)=(1-|z|)^\alpha$, $\alpha>-1$. Throughout this paper, we consider positive radial weights of the form $\omega=e^{-\eta}$ where $\eta$ is a strictly increasing radial function and $\Delta\eta>0$ on $\D$. Now, it is said that the weight $\omega$ belongs to the class $\mathcal{W}$ if we can choose a differentiable radial function $\tau$ as
\begin{align}\label{taudef}
(\Delta\eta(z))^{-\frac{1}{2}}\asymp\tau(r)
\end{align}
where $\lim_{r\rightarrow1^-}\tau(r)=0$ and $\lim_{r\rightarrow1^-}\tau'(r)=0$. We remark that $\mathcal{W}$ is a subclass of the so-called fast weights considered by Kriete and MacCluer in \cite{CM,KM}. Especially, $\eta(z)=\frac{1}{1-|z|}$ is contained in $\mathcal{W}$ as a typical example. Readers can refer to \cite{HLS,PP} to study other examples belonging to $\mathcal{W}$ and why these conditions are needed. On the other hand, we note that the standard weights $\eta(r)=-\alpha\log(1-r)$, $\alpha>-1$, are not contained in the class $\mathcal{W}$ since we can not find $\tau$ such that $\tau'(r)$ converge to $0$ when $r\rightarrow1^-$. That means we do not consider the standard weighted Bergman spaces in this paper.\\
\indent It is well known that all composition operators are bounded in the standard weighted Bergman spaces $A^p_\alpha(\D)$ by the Littlewood subordination principle, but it is not true in Bergman spaces with fast weights any more. Much studies on the theory of composition operators over the large Bergman spaces has been established in \cite{CM,KM}. When $\omega$ belongs to the class $\mathcal{W}$, it has known that $C_\varphi$ is bounded on $A^2(\omega)$  if and only if
\begin{align}\label{bddwithoutU}
\limsup_{|z|\rightarrow1}\frac{\omega(z)}{\omega(\varphi(z))}<\infty.
\end{align}
You can also refer to the proof in \cite[Theorem 3.2]{P} for the nessecity. As a consequence of the Carleson measure theorem, we observe that the condition (\ref{bddwithoutU}) still holds for the boundedness of $C_\varphi$ on $A^p(\omega)$ for all range of $p$ (See Remark \ref{bddp}). For the compactness, Kriete and MacCluer gave the estimate of the essential norm of $C_\varphi$ on $A^2(\omega)$:
\begin{align*}
\|C_\varphi\|_e\approx\limsup_{|z|\rightarrow1}\frac{\omega(z)}{\omega(\varphi(z))}.
\end{align*}

In the case of classical holomorphic function spaces such as the Hardy and the standard weighted Bergman space over the disk or the ball, many characterizations for the compactness of $C_\varphi-C_\psi$ have been developed over the past decades involving the pseudo-hyperbolic distance, but no result has been given in the large Bergman space; see for example \cite{CCKY, KW, M, MOZ}. In order to state our main result, we introduce the following distance induced by function $\tau$ associated with the weight $\omega\in\mathcal{W}$ as in (\ref{taudef}):
\begin{align}\label{ourdist}
\rho_{\tau,\varphi,\psi}(z):=\rho_\tau(\varphi(z),\psi(z))=1-e^{-d_\tau(\varphi(z),\psi(z))}
\end{align}
where
\begin{align}\label{dist}
d_\tau(z,w):=\inf_\gamma\int^1_0\frac{|\gamma'(t)|}{\tau(\gamma(t))}dt,
\end{align}
where the infimum is taken over all piecewise smooth curves $\gamma$ connecting $z$ and $w$. One may check that $\rho_{\tau}(z,w)$ can be a distance in \cite[Lemma 3.3]{P1}. In fact, the current author gave in \cite{P1} a sufficient condition for the compactness of $C_\varphi-C_\psi$ on $A^2(\omega)$ in terms of a distance induced by $\eta'$ function. In Section 3, we obtain the equivalent condition for compact difference of composition operators in terms of the distance involving $\tau$ function:
\begin{thm}\label{k-thm1}
Let $\omega\in\mathcal{W}$ and $C_\varphi$, $C_\psi$ be bounded on $A^q(\omega)$ for some $0<q<\infty$. Then the following statements are equivalent: for $0<p<\infty$,
\begin{itemize}
  \item [(1)]$C_\varphi-C_\psi$ is compact on $A^p(\omega)$.
  \item [(2)]$\lim_{|z|\rightarrow1^-}\rho_{\tau,\varphi,\psi}(z)\left(\frac{\omega(z)}{\omega(\varphi(z))}+\frac{\omega(z)}{\omega(\psi(z))}\right)=0$.
  \item [(3)]$\rho_{\tau,\varphi,\psi} C_\varphi$ and $\rho_{\tau,\varphi,\psi} C_\psi$ are compact from $A^p(\omega)$ into $L^p(\omega dA)$.
\end{itemize}
\end{thm}

One of the difficulties compared to the case of the standard weighted Bergman spaces is non-existence of the explicit form of the reproducing kernel in $A^2(\omega)$. In this paper, we overcome the difficulty using the Skwarczy\`nski distance which is defined by the reproducing kernel $\overline{K_z(w)}=K(z,w)$ in $A^2(\omega)$:
\begin{align*}
\mathcal{S}(z,w):=\sqrt{1-\frac{|K(z,w)|}{\|K_z\|\|K_w\|}},\quad z,w\in\D.
\end{align*}

We have proved that $\rho_\tau(z,w)$ is comparable to the distance $\mathcal{S}(z,w)$ (Theorem \ref{key1}) and the norm estimate for the difference of the reproducing kernels (Theorem \ref{diffkernelnorm}). More information on the Skwarczy\'nski distance will be presented in Section 2.2 and other essential materials to be used for this paper will be stated in Section 2. \\
\indent In Section 4, when the inducing maps have good boundary behavior in the sense of higher-order data and order of contact, we show that their difference is compact on $A^2(\omega)$. Using this result, we construct explicit analytic maps $\varphi, \psi$ which $C_\varphi$ and $C_\psi$ are non-compact composition operators but $C_\varphi-C_\psi$ is compact on $A^2(\omega)$. In Section 5, we characterize a topological path component of $\mathcal{C}_{HS}(A^2(\omega))$, the space of all composition operators on $A^2(\omega)$ endowed with topology induced by the Hilbert-Schmidt norm.
\begin{thm}\label{k-thm3}
Let $\omega\in\mathcal{W}$. Then the following statements are equivalent:
\begin{itemize}
  \item [(1)]$\rho_{\tau,\varphi,\psi}(z)\|K_{\varphi(z)}\|, \rho_{\tau,\varphi,\psi}(z)\|K_{\psi(z)}\|$ are belong to $L^2(\omega dA)$.
  \item [(2)]$C_\varphi-C_\psi$ is a Hilbert-Schmidt operator on $A^2(\omega)$.
  \item [(3)]$C_\varphi$ and $C_\psi$ lie in the same path component of $\mathcal{C}_{HS}(A^2(\omega))$.
\end{itemize}
\end{thm}

\bigskip

\textit{Constants.}
 In the rest of this paper, we use the notation
$X\lesssim Y$  or $Y\gtrsim X$ for nonnegative quantities $X$ and $Y$ to mean
$X\le CY$ for some inessential constant $C>0$. Similarly, we use the notation $X\approx Y$ if both $X\lesssim Y$ and $Y \lesssim X$ hold.

\bigskip
\section{Preliminary}
Assuming that $\tau'(r)\rightarrow0$ as $r\rightarrow1^-$ in Section 1, there exist constants $c_1, c_2>0$ independent of $z,w$ such that $\tau(z)< c_1(1-|z|)$ and
\begin{align}\label{comparable}
|\tau(z)-\tau(w)|\leq c_2|z-w|,\qquad z,w\in\D.
\end{align}
Throughout this paper, we denote
\begin{align*}
m_\tau:=\frac{\min(1,c_1^{-1},c_2^{-1})}{4}.
\end{align*}

\subsection{Radius functions and associated distances}
We let $D(z,r)$ be a Euclidean disk centred at $z$ with radius $r>0$. Using (\ref{comparable}), we obtain that for $0<\delta\leq m_\tau$,
\begin{align}\label{equiquan}
\frac{1}{2}\tau(z)\leq\tau(w)\leq2\tau(z)\quad\text{if}\quad w\in D(z,\delta\tau(z)),
\end{align}
and we use the notation $D(\delta\tau(z)):=D(z,\delta\tau(z))$ for simplicity.
One may refer to the proof of Lemma 2.1 of \cite{PP} for the inequality above. We set all $\delta$ appearing in the rest of our paper to meet the above conditions.

In \cite[page 12]{HLS}, when $\omega\in\mathcal{W}$, the authors gave the following useful inequality to estimate the reproducing kernel function of $A^2(\omega)$: for each positive number $M>0$, there exists a constant $C(M)>1$ such that
\begin{align}\label{estimate}
e^{-d_\tau(z,w)}\leq C\left(\frac{\min(\tau(z),\tau(w))}{|z-w|}\right)^M,\quad z\neq w\in\D.
\end{align}
Using (\ref{estimate}), we show that $d_\tau(z,w)$ and $\frac{|z-w|}{\tau(z)}$ are comparable whenever $d_\tau(z,w)<R$ for some $R>0$.
\begin{lem}\label{setinclus}
Let $\omega\in\mathcal{W}$. If $d_\tau(z,w)<R$ for some $R>0$ then there is a constant $C_1>0$ depending only on constants $M,C$ in (\ref{estimate}) and $R>0$ such that
\begin{align*}
d_\tau(z,w)\geq C_1\frac{|z-w|}{\min(\tau(z),\tau(w))}.
\end{align*}
\end{lem}
\begin{proof}
For a given $z\in\D$ and any point $w\neq z$ with $d_\tau(z,w)<R$, from (\ref{estimate})
\begin{align*}
\frac{|z-w|}{\min(\tau(z),\tau(w))}<(Ce^R)^{\frac{1}{M}}.
\end{align*}
Thus, if we assume that $|z|\geq|w|$, then we have a constant $0<s\leq1/\delta(Ce^R)^{\frac{1}{M}}$ satisfying
\begin{align}\label{distzw}
|z-w|=s\delta\tau(z).
\end{align}
Letting $\gamma$ denote the curve obtained so that $\gamma(0)=z$, $\gamma(1)=w$ and $d_\tau(z,w)=\int^{1}_0\frac{|\gamma'(t)|}{\tau(\gamma(t))}dt$, we choose the minimum value $0<t_0\leq1$ such that
\begin{align}\label{distzw1}
|z-\gamma(t_0)|=\delta\tau(z).
\end{align}
For $1<s\leq1/\delta(Ce^R)^{\frac{1}{M}}$, by (\ref{equiquan}), (\ref{distzw}) and (\ref{distzw1}), we have
\begin{align}\label{s<1}
d_\tau(z,w)\geq\int^{t_0}_0\frac{|\gamma'(t)|}{\tau(\gamma(t))}dt\geq\frac{1}{2\tau(z)}\int^{t_0}_0|\gamma'(t)|dt&\geq\frac{\delta}{2}=\frac{|z-w|}{2s\tau(z)}\\
&\geq\frac{\delta}{2}(Ce^R)^{-\frac{1}{M}}\frac{|z-w|}{\tau(z)}.\nonumber
\end{align}
For $0<s\leq1$, we obtain $d_\tau(z,w)\geq\frac{|z-w|}{2\tau(z)}$ promptly from (\ref{s<1}). Thus, if we choose $C_1=\delta/2(Ce^R)^{-\frac{1}{M}}$ then we obtain the desired inequality.
\end{proof}

\subsection{Skwarczy\`nski distance $\mathcal{S}(z,w)$}
Given $\omega\in\mathcal{W}$, the space $A^2(\omega)$ is a closed subspace of $L^2(\omega dA)$ with inner product
\begin{align*}
\langle f,g\rangle_\omega=\int_\D f\overline{g}\omega^2 dA,\quad f,g\in A^2(\omega).
\end{align*}
As is well known, the reproducing kernel for the Bergman space $A^2(\omega)$ is defined by
\begin{align}\label{definition}
K(z,w)=\overline{K_z(w)}=\sum_{k=0}^\infty e_k(z)\overline{e_k(w)}
\end{align}
where $\{e_k\}$ is an arbitrary orthonormal basis for $A^2(\omega)$. The reproducing kernel estimates of Bergman spaces with exponential weights $\omega=e^{-\eta}$ have been established by \cite{AH,AP,HLS,LR} for example: for $\omega\in\mathcal{W}$ and $\tau(z)\approx(\Delta\eta(z))^{-1/2}$, there are positive constants $\sigma, C>0$ such that
\begin{align}\label{kernelesti}
\begin{aligned}
&|K(z,w)|\omega(z)\omega(w)\leq \frac{C}{\tau(z)\tau(w)}e^{-\sigma d_\tau(z,w)},\quad\forall z,w\in\D,\\
&\int_\D|K(z,\cdot)\omega|^pdA\asymp\omega(z)^{-p}\tau(z)^{-2(p-1)},\quad0<p<\infty.
\end{aligned}
\end{align}
 Now, we introduce another distance defined using the reproducing kernel function:
\begin{align*}
\mathcal{S}(z,w):=\left(1-\frac{|K(z,w)|}{\|K_z\|\|K_w\|}\right)^{1/2},\quad z,w\in\D.
\end{align*}
It is called the Skwarczy\`nski pseudodistance and we can refer to \cite{JP} for studying about the pseudodistance $\mathcal{S}(z,w)$, especially its relation with the Bergman pseudodistance $\beta(z,w)$ induced by the Bergman pseudometric,
\begin{align*}
B\left(\frac{\partial}{\partial z},z\right)=\frac{\sup\{|f'(z)|:f\in A^2(\omega), f(z)=0, \|f\|=1\}}{\sqrt{K(z,z)}};
\end{align*}
for example, we have the following inequality in \cite[Corollary 6.4.7 (a)]{JP}:
\begin{align}\label{relbergman}
\mathcal{S}(z,w)\leq \frac{1}{\sqrt{2}}\beta(z,w),\quad z,w\in\D.
\end{align}
Moreover, we have the following inequality in \cite[Theorem 6.4.3]{JP} to compare with our distance $\rho_\tau(z,w)$:
\begin{align}\label{lowerds}
\mathcal{S}(z,w)\leq\frac{M(z,w)}{\sqrt{K(z,z)}}\leq\sqrt{2}\mathcal{S}(z,w),\quad z,w\in\D
\end{align}
where $M(z,w):=\sup\{|f(z)|:f\in A^2(\omega), \|f\|=1, f(w)=0\}$.
\begin{lem}\label{caltestft}
Given $z, w\in\D$ with $d_\tau(z,w)<R$ for some $R>0$, define
\begin{align*}
f_{z,w}(\xi)=\omega(z)K(\xi,z)(\xi-w)
\end{align*}
then $\|f_{z,w}\|\leq C_2$ where a constant $C_2:=C_2(R)>0$.
\end{lem}
\begin{proof}
Given $0<\delta\leq m_\tau$, we decompose the disk into two parts as follows:
\begin{align}\label{goal5}
\|f_{z,w}\|^2&=\omega(z)^2\int_\D|K(\xi,z)|^2|\xi-w|^2\omega(\xi)^2dA(\xi)\nonumber\\
&=\int_{D(\delta\tau(z))}+\int_{\D\setminus D(\delta\tau(z))}\omega(z)^2|K(\xi,z)|^2|\xi-w|^2\omega(\xi)^2dA(\xi).
\end{align}
Since $|z-w|<R/C_1\tau(z)$ where $C_1$ is defined in Lemma \ref{setinclus}, (\ref{kernelesti}) and the triangle inequality yield that
\begin{align*}
&\omega(z)^2\int_{D(\delta\tau(z))}|K(\xi,z)|^2|\xi-w|^2\omega(\xi)^2dA(\xi)\\
&\lesssim\frac{1}{\tau(z)^2}\int_{D(\delta\tau(z))}\frac{|\xi-z|^2+|z-w|^2}{\tau(\xi)^2}dA(\xi)\lesssim1.
\end{align*}
On the other hand, by Lemma \ref{setinclus} with $M=3/\sigma$, we have
\begin{align*}
|z-w|<\frac{R}{C_1}\tau(z)<2R(Ce^R)^{\frac{\sigma}{3}}|z-\xi|\quad\text{for}\quad\xi\in\D\setminus D(\delta\tau(z)).
\end{align*}
This, together with (\ref{estimate}) and (\ref{kernelesti}), gives the following inequality:
\begin{align*}
&\omega(z)^2\int_{\D\setminus D(\delta\tau(z))}|K(\xi,z)|^2|\xi-w|^2\omega(\xi)^2dA(\xi)\\
&\lesssim\int_{\D\setminus D(\delta\tau(z))}\frac{|\xi-z|^2+|z-w|^2}{\tau(z)^2\tau(\xi)^2}\left(\frac{\min(\tau(z),\tau(\xi))}{|\xi-z|}\right)^6dA(\xi)\\
&\lesssim\int_{\D\setminus D(\delta\tau(z))}\frac{\tau(z)^2}{|\xi-z|^4}dA(\xi)\leq\sum_{j=0}^{\infty}\int_{2^j\delta\tau(z)<|z-\xi|\leq2^{j+1}\delta\tau(z)}\frac{\tau(z)^2}{|\xi-z|^4}dA(\xi)\\
&\lesssim\sum_{j=0}^{\infty}\frac{1}{2^{4j}\delta^4\tau(z)^2}2^{2j+2}\delta^2\tau(z)^2\lesssim1.
\end{align*}
\end{proof}
The comparison of the following two distances is essential to complete our main results, Theorem \ref{k-thm1} and Theorem \ref{k-thm3}.

\begin{thm}\label{key1}
Let $\omega\in\mathcal{W}$. Then $\rho_\tau(z,w)\approx\mathcal{S}(z,w)$ for $z,w\in\D$.
\end{thm}
\begin{proof}
By the kernel estimate (\ref{kernelesti}), there is a constant $C_3>0$ such that
\begin{align*}
\frac{|K(z,w)|}{\|K_z\|\|K_w\|}\leq C_3e^{-\sigma d_\tau(z,w)},\quad z,w\in\D.
\end{align*}
From the inequality above, when $d_\tau(z,w)\geq R$ with $C_3e^{-\sigma R}<1/4$,
\begin{align*}
\frac{1}{2}\leq\frac{(1-C_3e^{-\sigma d_\tau(z,w)})^{1/2}}{1-e^{-d_\tau(z,w)}}\leq\frac{\mathcal{S}(z,w)}{\rho_\tau(z,w)}.
\end{align*}
Now, let's consider the case of $d_\tau(z,w)<R$ with $\tau(z)\leq\tau(w)$. Applying $f_{z,w}(\xi)=\omega(z)K(\xi,z)(\xi-w)$ to (\ref{lowerds}), we have
\begin{align*}
\mathcal{S}(z,w)\geq\frac{1}{\sqrt{2}}\frac{|f_{z,w}(z)|}{\sqrt{K(z,z)}}\gtrsim\frac{|z-w|}{\tau(z)}\geq d_\tau(z,w)\geq\rho_\tau(z,w).
\end{align*}
Therefore, $\rho_\tau(z,w)\lesssim\mathcal{S}(z,w)$ for any $z,w\in\D$. Meanwhile, since it is known that $d_\tau$ is comparable to the Bergman distance $\beta$ by \cite[p.355]{CO}, together with (\ref{relbergman}), we have a constant $C>1$ independent of $z,w$ such that
\begin{align*}
\mathcal{S}(z,w)\leq e\big[1-e^{-\mathcal{S}(z,w)}\big]&\lesssim1-e^{-1/2\beta(z,w)}\\&\leq1-e^{-Cd_\tau(z,w)}\leq C\rho_\tau(z,w).
\end{align*}
Therefore, we complete our proof.
\end{proof}

\subsection{Carleson measure theorem}
A positive Borel measure $\mu$ in $\D$ is called a (vanishing) Carleson measure for $A^p(\omega)$ if the embedding $A^p(\omega)\subset L^p(\omega d\mu)$ is (compact) continuous where
\begin{align}\label{embedding}
L^p(\omega d\mu):=\left\{f\in\mathcal{M}(\D)\big|\int_\D|f(z)\omega(z)|^pd\mu(z)<\infty\right\}
\end{align}
and $\mathcal{M}(\D)$ is a set of $\mu$-measurable functions on $\D$. Now, we introduce Carleson measure theorem in our setting which can be found in \cite{LR, O, PP} for example.
\begin{thm}[Carleson measure theorem]
Let $\omega\in\mathcal{W}$ and $\mu$ be a positive Borel measure on $\D$. Then, for $0<p<\infty$, we have
\begin{itemize}
\item[(1)]The embedding $I:A^p(\omega)\rightarrow L^p(\omega d\mu)$ is bounded if and only if for a small $\delta\in(0,m_\tau)$, we have $\sup_{z\in\D}\frac{\mu(D(\delta\tau(z)))}{\tau(z)^{2}}<\infty$.
\item[(2)]The embedding $I:A^p(\omega)\rightarrow L^p(\omega d\mu)$ is compact if and only if for a small $\delta\in(0,m_\tau)$, we have $\lim_{|z|\rightarrow1}\frac{\mu(D(\delta\tau(z)))}{\tau(z)^{2}}=0$.
\end{itemize}
\end{thm}
By the measure theoretic change of variables, we have
\begin{align*}
\|uC_\varphi f\|_{p}^p=\int_\D|u(f\circ\varphi)\omega|^p dA=\int_\D|f\omega|^p d\mu_{u,\varphi,p}
\end{align*}
where
\begin{align*}
\mu_{u,\varphi,p}(E):=\omega^{-p}[|u|^p\omega^{p}dA]\circ\varphi^{-1}(E)=\int_{\varphi^{-1}(E)}|u(z)|^p\frac{\omega(z)^{p}}{\omega(\varphi(z))^{p}}dA(z)
\end{align*}
for any measurable subsets $E$ of $\D$. Denote the norm of operator $T$ acting on $A^p(\omega)$ by $\|T\|_{A^p(\omega)}$. It is well known that
\begin{align}\label{opnorm}
\|uC_\varphi\|^p_{A^p(\omega)}\approx\sup_{z\in\D}\frac{\mu_{u,\varphi,p}(D(\delta\tau(z)))}{\tau(z)^{2}}.
\end{align}
Here, we can obtain the following relation between $\|C_\varphi\|_{A^2(\omega)}$ and the norm of $C_\varphi$ acting on $A^p(\omega^{2/p})$:
\begin{align}\label{bddindp}
\|C_\varphi\|^p_{A^p(\omega^{2/p})}\approx\sup_{z\in\D}\frac{\omega^{-2}[\omega^2dA]\circ\varphi^{-1}(D(\delta\tau(z)))}{\tau(z)^{2}}\approx\|C_\varphi\|^2_{A^2(\omega)}
\end{align}
This holds because the weight $\omega^{2/p}$ still belongs to $\mathcal{W}$, and we can choose the same $\tau$ functions for $A^p(\omega^{2/p})$ with those for $A^2(\omega)$. As a consequence of (\ref{bddindp}), we conclude the following Remark:
\begin{rmk}\label{bddp}
Let $\omega\in\mathcal{W}$, $0<p<\infty$. Then $C_\varphi$ is bounded on $A^p(\omega)$ if and only if
\begin{align*}
\sup_{z\in\D}\frac{\omega(z)}{\omega(\varphi(z))}<\infty.
\end{align*}
Accordingly, letting $\varphi_s(z)=(1-s)\varphi(z)+s\psi(z)$ for $0\leq s\leq1$, $C_{\varphi_s}$ is uniformly bounded on $A^p(\omega)$ with respect to $s$ when $C_\varphi$, $C_\psi$ are bounded on $A^p(\omega)$.
\end{rmk}
\begin{proof}
From (\ref{bddindp}), the boundedness of $C_\varphi$ acting on $A^p(\omega)$ is equivalent to the boundedness of $C_\varphi$ acting on $A^2(\omega^{p/2})$. Thus, it is clear that the first assertion holds by (\ref{bddwithoutU}). Furthermore, since $|\varphi_s(z)|\leq\max(|\varphi(z)|,|\psi(z)|)$, we have
\begin{align}\label{uniformbdd}
\frac{\omega(z)}{\omega(\varphi_s(z))}\leq\frac{\omega(z)}{\omega(\varphi(z))}+\frac{\omega(z)}{\omega(\psi(z))}.
\end{align}
Therefore, $C_{\varphi_s}$ is uniformly bounded with respect to $0\leq s\leq1$ when $C_\varphi$, $C_\psi$ are bounded on $A^p(\omega)$ by (\ref{bddwithoutU}).
\end{proof}

\section{Difference of Composition operators}

In this section, we prove Theorem \ref{k-thm1} and other interests. The following lemma which gives the upper estimate for the derivative of $f$ plays a crucial role in our proofs. In fact, its proof is similar to the case of a doubling measure $\Delta\eta$ whose proof we can find in \cite[Lemma 19]{MMO}. For the proofs of our setting, you can refer to \cite{HLS,O,PP}.
\begin{lem}\label{derivsubmean}
Let $\omega=e^{-\eta}$, where $\eta$ is a subharmonic function and $0<p<\infty$. Suppose the function $\tau$ satisfies properties (\ref{equiquan}) and  $\tau(z)^2\Delta\eta(z)\lesssim1$. Given $\delta>0$ satisfying (\ref{equiquan}) and $f\in H(\D)$,
\begin{itemize}
  \item[(1)]$|f(z)|^pe^{-\beta\eta(z)}\lesssim\frac{1}{\tau(z)^2}\int_{D(\delta\tau(z))}|f(\xi)|^pe^{-\beta\eta(\xi)}dA(\xi)$,
  \item[(2)]$|f'(z)e^{-\eta(z)}|^p\lesssim\frac{1}{\tau(z)^{2+p}}\int_{D(\delta\tau(z))}|f(\xi) e^{-\eta(\xi)}|^p\,dA(\xi)$.
\end{itemize}
\end{lem}
First, we prove the following inequality which gives the upper estimate for the difference in function values for two variables $a, b$ with $|a-b|<\delta\min(\tau(a),\tau(b))$.
\begin{lem}\label{mvtforinterg}
Let $\omega=e^{-\eta}\in\mathcal{W}$ and $0<p<\infty$. Then there is a constant $C:=C(\delta,p)>0$ such that
\begin{align*}
|f(a)-f(b)|^pe^{-p\eta(a)}\leq C\frac{\rho_\tau(a,b)^p}{\tau(a)^{2}}\int_{D(6\delta\tau(a))}|f(\xi) e^{-\eta(\xi)}|^p\,dA(\xi),
\end{align*}
for all $f\in H(\D)$ and all $b\in D(\delta\tau(a))$ with $|a|\geq|b|$.
\end{lem}
\begin{proof}
By the fundamental theorem of calculus, for $b\in D(\delta\tau(a))$,
\begin{align*}
|f(a)-f(b)|&=\left|\int_{0}^{1}f'(at+b(1-t))(a-b)dt\right|\\&\leq|a-b|\sup_{t\in[0,1]}|f'(at+b(1-t))|.
\end{align*}
Denote $a_t:=at+b(1-t)$, $0\leq t\leq1$. By (\ref{equiquan}), for $z\in D(\delta\tau(a_t))$, we have
\begin{align}\label{setinclus1}
|z-a_i|&\leq|z-a_t|+|a_t-a_i|\nonumber\\&\leq\delta\tau(a_t)+\delta\tau(a)\leq3\delta\tau(a)\leq6\delta\tau(a_i)\quad\text{for}\quad i=0,1.
\end{align}
Since $|a_t|\leq\max(|a|,|b|)=|a|$ by our assumption, (2) of Lemma \ref{derivsubmean} and (\ref{setinclus1}) yield that
\begin{align*}
|f'(a_t)|^p&\lesssim\frac{e^{p\eta(a_t)}}{\tau(a_t)^{2+p}}\int_{D(\delta\tau(a_t))}|f(\xi) e^{-\eta(\xi)}|^p\,dA(\xi)\\
&\lesssim\frac{e^{p\eta(a)}}{\tau(a)^{2+p}}\int_{D(6\delta\tau(a))}|f(\xi) e^{-\eta(\xi)}|^p\,dA(\xi).
\end{align*}
Thus, together with Lemma \ref{setinclus} and $d_\tau(a,b)<2\delta$, there is a constant $C:=C(\delta,p,M)>0$ such that
\begin{align*}
|f(a)-f(b)|^pe^{-p\eta(a)}\leq C\frac{d_\tau(a,b)^p}{\tau(a)^{2}}\int_{D(6\delta\tau(a))}|f(\xi) e^{-\eta(\xi)}|^p\,dA(\xi).
\end{align*}
Since $e^{-2\delta}x\leq1-e^{-x}$ for $0\leq x<2\delta<<1$, we have
\begin{align*}
d_\tau(a,b)\leq e^{2\delta}\rho_{\tau}(a,b).
\end{align*}
Hence, we complete the proof.
\end{proof}
Using the result above, we get the following optimal norm estimate for the difference of the reproducing kernels:
\begin{thm}\label{diffkernelnorm}
Let $\omega\in\mathcal{W}$. Then for any $z,w\in\D$,
\begin{align*}
\frac{\|K_z-K_w\|^2}{\|K_z\|^2+\|K_w\|^2}\approx\rho_\tau(z,w)^2.
\end{align*}
\end{thm}
\begin{proof}
By Theorem \ref{key1}, we have
\begin{align*}
\|K_{z}-K_{w}\|^2&\geq\|K_{z}\|^2+\|K_{w}\|^2-2|K(z,w)|\\
&=\|K_{z}\|^2+\|K_{w}\|^2-2\|K_{z}\|\|K_{w}\|(1-\mathcal{S}(z,w)^2)\\
&\geq(\|K_{z}\|^2+\|K_{w}\|^2)\mathcal{S}(z,w)^2\\
&\gtrsim(\|K_{z}\|^2+\|K_{w}\|^2)\rho_\tau(z,w)^2.
\end{align*}
Conversely, we first consider the case $d_\tau(z,w)<R$ where $0<R<\delta^2/12(Ce)^{-\frac{1}{M}}$, $C, M$ appeared in (\ref{estimate}). Then $|z-w|<\delta/6\tau(z)$ by Lemma \ref{setinclus}, so for $|z|\geq|w|$,
\begin{align*}
\|K_{z}-K_{w}\|^2&=\sum_{n=0}^{\infty}\frac{1}{\|z^n\|^2}|z^n-w^n|^2\nonumber\\
&\lesssim\frac{\rho_{\tau}(z,w)^2e^{2\eta(z)}}{\tau(z)^2}
\sum_{n=0}^{\infty}\frac{1}{\|z^n\|^2}\int_{D(\delta\tau(z))}|\xi|^{2n}e^{-2\eta(\xi)}dA(\xi),
\end{align*}
by Proposition \ref{mvtforinterg}. Meanwhile, (\ref{definition}), (\ref{kernelesti}) and (\ref{equiquan}) yield that
\begin{align*}
\sum_{n=0}^{\infty}\frac{1}{\|z^n\|^2}\int_{D(\delta\tau(z))}|\xi|^{2n}e^{-2\eta(\xi)}dA(\xi)
&=\int_{D(\delta\tau(z))}\|K_\xi\|^2e^{-2\eta(\xi)}dA(\xi)\\
&\approx\int_{D(\delta\tau(z))}\frac{1}{\tau(\xi)^2}dA(\xi)\approx1.
\end{align*}
This, together with the kernel estimate (\ref{kernelesti}), yields that
\begin{align*}
\|K_{z}-K_{w}\|^2\lesssim\rho_{\tau}(z,w)^2\|K_z\|^2,\quad\text{for}\quad|z|\geq|w|.
\end{align*}
In the same way, we have
\begin{align*}
\|K_{z}-K_{w}\|^2\lesssim\rho_{\tau}(z,w)^2\|K_w\|^2,\quad\text{for}\quad|z|<|w|.
\end{align*}
For the case $d_\tau(z,w)\geq R$, we easily obtain from the triangle inequality that
\begin{align*}
\|K_{z}-K_{w}\|^2\lesssim\rho_{\tau}(z,w)^2(\|K_z\|^2+\|K_w\|^2).
\end{align*}
Therefore, we complete our proof.
\end{proof}
As mentioned earlier, not every $C_\varphi$ is bounded on $A^p(\omega)$, so we can give the question when the difference of two composition operators is bounded on $A^p(\omega)$. In Proposition \ref{bddnecess} and Proposition \ref{mainresult2}, we give  characterizations for the boundedness of $C_\varphi-C_\psi$ on $A^p(\omega)$.
\begin{lem}\cite[Lemma 4.1]{P1}\label{delight}
Let $\omega\in\mathcal{W}$ and $\varphi\in S(\D)$. If there is a curve $\gamma$ connecting to $\zeta\in\partial\D$ and a constant $c>0$ such that $\liminf_{z\to \zeta}\frac{\omega(z)}{\omega(\varphi(z))}\geq c$ where $z\in\gamma$, then
\begin{align*}
\liminf_{z\to \zeta}\frac{\tau(z)}{\tau(\varphi(z))}\geq\min(1,c)\quad\text{for}\quad z\in\gamma.
\end{align*}
Moreover, if $\lim_{z\to \zeta}\frac{\omega(z)}{\omega(\varphi(z))}=0$ then $\limsup_{z\to \zeta}\frac{\tau(z)}{\tau(\varphi(z))}\leq1$.
\end{lem}

\begin{pro}\label{bddnecess}
Let $\omega\in\mathcal{W}$. If $C_\varphi-C_\psi$ is bounded on $A^2(\omega)$, then
\begin{align*}
\sup_{z\in\D}\rho_{\tau,\varphi,\psi}(z)\left(\frac{\omega(z)}{\omega(\varphi(z))}+\frac{\omega(z)}{\omega(\psi(z))}\right)<\infty.
\end{align*}
\end{pro}
\begin{proof}
Assume that there exists a sequence $\{z_n\}$ such that
\begin{align}\label{assump1}
\rho_{\tau,\varphi,\psi}(z_n)\frac{\omega(z_n)}{\omega(\varphi(z_n))}\rightarrow\infty\quad\text{as}\quad n\rightarrow\infty.
\end{align}
Then we easily see that $\lim_{n\rightarrow\infty}\frac{\omega(z_n)}{\omega(\varphi(z_n))}=\infty$. By Theorem \ref{diffkernelnorm}, Lemma \ref{delight} and the kernel estimate,
\begin{align*}
\infty>\|C_\varphi^*-C_\psi^*\|_{A^2(\omega)}^2&\geq\frac{1}{\|K_{z_n}\|^2}\|K_{\varphi(z_n)}-K_{\psi(z_n)}\|^2\\&\gtrsim\rho_{\tau,\varphi,\psi}(z_n)^2\left(\frac{\|K_{\varphi(z_n)}\|^2}
{\|K_{z_n}\|^2}+\frac{\|K_{\psi(z_n)}\|^2}{\|K_{z_n}\|^2}\right)\\
&\gtrsim\rho_{\tau,\varphi,\psi}(z_n)^2\frac{\omega(z_n)^2}
{\omega(\varphi(z_n))^2}.
\end{align*}
This yields a contradiction to (\ref{assump1}) when $n\rightarrow\infty$.
\end{proof}

In what follows, we denote the operator norm $\|T\|_{A^p(\omega)}$ by $\|T\|_p$ for simplicity.
\begin{pro}\label{mainresult2}
Let $\omega=e^{-\eta}\in\mathcal{W}$. Then, for $0<p<\infty$,
\begin{align*}
\|C_\varphi-C_\psi\|_{p}\lesssim\|\rho_{\tau,\varphi,\psi} C_\varphi\|_{p}+\|\rho_{\tau,\varphi,\psi} C_\psi\|_{p}.
\end{align*}
\end{pro}
\begin{proof}
For $\|f\|_p\leq1$, we write
\begin{align*}
&\|(C_\varphi-C_\psi)f\|^p_p\\
&=\int_{E}|(f\circ\varphi-f\circ\psi) e^{-\eta}|^pdA+\int_{E^c}|(f\circ\varphi-f\circ\psi)e^{-\eta}|^pdA
\end{align*}
where
\begin{align}\label{set}
E:=\{z\in\D:d_\tau(\varphi(z),\psi(z))<R<1\}.
\end{align}
Take $0<R<\delta^2/12(Ce)^{-\frac{1}{M}}$ where $M,C$ in (\ref{estimate}). Then by Lemma \ref{setinclus},
\begin{align*}
|\varphi(z)-\psi(z)|<\delta/6\tau(\varphi(z))\quad\text{for}\quad z\in E.
\end{align*}
Also, we denote $E_1:=\{z\in E:|\varphi(z)|\geq|\psi(z)|\}$. Then by Lemma \ref{mvtforinterg}, we have
\begin{align*}
&\int_{E}|(f\circ\varphi-f\circ\psi)e^{-\eta}|^pdA\\
&\lesssim\int_{E_1}\frac{\rho_{\tau,\varphi,\psi}(z)^pe^{p\eta(\varphi(z))-p\eta(z)}}{\tau(\varphi(z))^2}
\int_{D(\delta\tau(\varphi(z)))}|fe^{-\eta}|^pdA\,dA(z)\\
&+\int_{E\setminus E_1}\frac{\rho_{\tau,\varphi,\psi}(z)^pe^{p\eta(\psi(z))-p\eta(z)}}{\tau(\psi(z))^2}
\int_{D(\delta\tau(\psi(z)))}|fe^{-\eta}|^pdA\,dA(z).
\end{align*}

Therefore, by (\ref{equiquan}), (\ref{opnorm}) and Fubini's theorem, we have
\begin{align}
&\int_{E}|(f\circ\phi-f\circ\psi)e^{-\eta}|^pdA\nonumber\\
&\lesssim\int_{\D}|f(\xi)e^{-\eta(\xi)}|^p\int_{\varphi^{-1}(D(\delta\tau(\xi)))}
\frac{\chi_{E_1}(z)\rho_{\tau,\varphi,\psi}(z)^p}{\tau(\varphi(z))^2}\frac{e^{-p\eta(z)}}{e^{-p\eta(\varphi(z))}}dA(z)dA(\xi)\nonumber\\
&+\int_{\D}|f(\xi)e^{-\eta(\xi)}|^p\int_{\psi^{-1}(D(\delta\tau(\xi)))}
\frac{\chi_{E\setminus E_1}(z)\rho_{\tau,\varphi,\psi}(z)^p}{\tau(\psi(z))^2}\frac{e^{-p\eta(z)}}{e^{-p\eta(\psi(z))}}dA(z)dA(\xi)\nonumber\\
&\lesssim\|\chi_{E}\rho_{\tau,\varphi,\psi} C_\varphi\|^p_{p}+\|\chi_{E}\rho_{\tau,\varphi,\psi} C_\psi\|^p_{p}.\label{integralE}
\end{align}
Since $\rho_{\tau,\varphi,\psi}(z)\geq1-e^{-R}$ on $E^c$, we have
\begin{align}\label{integralEc}
\int_{E^c}|f\circ\varphi-f\circ\psi|^pe^{-p\eta} dA&\lesssim\int_{E^c}(|f\circ\varphi|^p+|f\circ\psi|^p)\rho_{\tau,\varphi,\psi}^pe^{-p\eta}dA\nonumber\\
&\leq\|\chi_{E^c}\rho_{\tau,\varphi,\psi} C_\varphi f\|^p_p+\|\chi_{E^c}\rho_{\tau,\varphi,\psi} C_\psi f\|^p_p.
\end{align}
By (\ref{integralE}) and (\ref{integralEc}), we complete our proof.
\end{proof}
\begin{rmk}
Combining Proposition \ref{bddnecess} with Proposition \ref{mainresult2}, we remark that if $\rho_{\tau,\varphi,\psi} C_\varphi$ is bounded from $A^2(\omega)$ into $L^2(\omega dA)$ then
\begin{align*}
\sup_{z\in\D}\rho_{\tau,\varphi,\psi}(z)\frac{\omega(z)}{\omega(\varphi(z))}<\infty.
\end{align*}
\end{rmk}
In the same way as the proof of Proposition \ref{mainresult2}, we can obtain the following sufficient condition for the compactness of $C_\varphi-C_\psi$ immediately.
\begin{pro}\label{mainresult3}
Let $\omega=e^{-\eta}\in\mathcal{W}$ and $0<p<\infty$. If $\rho_{\tau,\varphi,\psi} C_\varphi$ and $\rho_{\tau,\varphi,\psi} C_\psi$ are compact from $A^p(\omega)$ into $L^p(\omega dA)$ then $C_\varphi-C_\psi$ is compact on $A^p(\omega)$.
\end{pro}
\begin{proof}
Consider a sequence $\{f_k\}$ converging to $0$ weakly on $A^p(\omega)$ when $k\rightarrow\infty$. We claim that the following integral vanishes as $k\rightarrow\infty$,
\begin{align}\label{goal}
&\|f_k\circ\varphi-f_k\circ\psi\|_p^p\nonumber\\
&=\int_{E}|f_k\circ\varphi-f_k\circ\psi|^pe^{-p\eta}dA+\int_{E^c}|f_k\circ\varphi-f_k\circ\psi|^pe^{-p\eta}dA
\end{align}
where the set $E$ is defined in (\ref{set}). First, we easily see that the second integral of (\ref{goal}) vanishes as $k\rightarrow\infty$ by (\ref{integralEc}) of Proposition \ref{mainresult2} and the criterion for the compactness. Thus, we only need to verify that the first integral of (\ref{goal}) converges to $0$. Since $\rho_{\tau,\varphi,\psi} C_\varphi$ and $\rho_{\tau,\varphi,\psi} C_\psi$ are compact on $A^p(\omega)$, there is $r(\epsilon)>0$ such that for $r<|\xi|<1$,
\begin{align*}
\frac{1}{\tau(\xi)^2}\left(\int_{\varphi^{-1}(D(\delta\tau(\xi)))}
\rho_{\tau,\varphi,\psi}^p\frac{e^{-p\eta}}{e^{-p\eta\circ\varphi}}dA+\int_{\psi^{-1}(D(\delta\tau(\xi)))}
\rho_{\tau,\varphi,\psi}^p\frac{e^{-p\eta}}{e^{-p\eta\circ\psi}}dA\right)<\epsilon.
\end{align*}
Using the same method in Proposition \ref{mainresult2}, the first integral of (\ref{goal}) is dominated by
\begin{align}\label{integralEK}
&\int_{E}|(f_k\circ\varphi-f_k\circ\psi)e^{-\eta}|^pdA\nonumber\\
&\lesssim\int_{\D}|f_k(\xi)|^pe^{-p\eta(\xi)}\int_{\varphi^{-1}(D(\delta\tau(\xi)))}
\frac{\chi_{E_1}(z)\rho_{\tau,\varphi,\psi}(z)^p}{\tau(\varphi(z))^2}\frac{e^{-p\eta(z)}}{e^{-p\eta(\varphi(z))}}dA(z)dA(\xi)\nonumber\\
&+\int_{\D}|f_k(\xi)e^{-\eta(\xi)}|^p\int_{\psi^{-1}(D(\delta\tau(\xi)))}
\frac{\chi_{E\setminus E_1}(z)\rho_{\tau,\varphi,\psi}(z)^p}{\tau(\psi(z))^2}\frac{e^{-p\eta(z)}}{e^{-p\eta(\psi(z))}}dA(z)dA(\xi)\nonumber\\
&\lesssim\epsilon\int_{\D\setminus r\D}|f_k|^pe^{-p\eta}dA+(\|\rho_{\tau,\varphi,\psi} C_\varphi\|_p^p+\|\rho_{\tau,\varphi,\psi} C_\psi\|_p^p)\int_{r\D}|f_k|^pe^{-p\eta}dA.
\end{align}
Therefore, we can make the integration above small when $k\rightarrow\infty$, so we complete our proof.
\end{proof}
We have shown $(3)\Longrightarrow(1)$ of Theorem \ref{k-thm1} in Proposition \ref{mainresult3}. Now, we prove the implication $(1)\Longrightarrow(2)$ of Theorem \ref{k-thm1}.
\begin{pro}\label{mainresult4}
Let $\omega\in\mathcal{W}$ and $C_\varphi$, $C_\psi$ be bounded on $A^p(\omega)$ for some $p\in(0,\infty)$. If $C_\varphi-C_\psi$ is compact on $A^p(\omega)$, then
\begin{align*}
\lim_{|z|\rightarrow1^-}\rho_{\tau,\varphi,\psi}(z)\left(\frac{\omega(z)}{\omega(\varphi(z))}+\frac{\omega(z)}{\omega(\psi(z))}\right)=0.
\end{align*}
\end{pro}
\begin{proof}
Assume that there is a boundary point $\zeta$ such that
\begin{align*}
\lim_{z\rightarrow\zeta}\rho_{\tau,\varphi,\psi}(z)\frac{\omega(z)}{\omega(\varphi(z))}\neq0.
\end{align*}
Then we can choose a sequence $\{z_n\}$ satisfying $\lim_{n\rightarrow\infty}\frac{\omega(z_n)}{\omega(\varphi(z_n))}\gtrsim1$ with $|\varphi(z_n)|\geq|\psi(z_n)|$ and $\lim_{n\rightarrow\infty}\rho_{\tau,\varphi,\psi}(z_n)=c>0$ from the boundedness of $C_\varphi$ and $C_\psi$. Now, define the bounded sequence $\{g_n\}$:
\begin{align}\label{testft}
g_n(\xi)=\frac{K_{\varphi(z_n)}(\xi)}{\omega(\varphi(z_n))^{-1}\tau(\varphi(z_n))^{-2+\frac{2}{p}}}
\end{align}
which uniformly converges to $0$ on compact subsets of $\D$. By (1) of Lemma \ref{derivsubmean}, Lemma \ref{delight} and Theorem \ref{key1},
\begin{align}\label{goal3}
\|(C_\varphi-C_{\psi})g_{n}\|_p^p&\geq\int_{D(\delta\tau(z_n))}|g_{n}(\varphi(\xi))-g_{n}(\psi(\xi))|^p\omega(\xi)^pdA(\xi)\nonumber\\
&\geq\tau(z_n)^2|g_n(\varphi(z_n))-g_{n}(\psi(z_n))|^p\omega(z_n)^p\nonumber\\
&\gtrsim|K_{\varphi(z_n)}(\varphi(z_n))-K_{\varphi(z_n)}(\psi(z_n))|^p\frac{\omega(z)^p}{\omega(\varphi(z_n))^{-p}}\frac{\tau(z_n)^2}{\tau(\varphi(z_n))^{2(1-p)}}\nonumber\\
&\gtrsim\frac{\omega(z_n)^p}{\omega(\varphi(z_n))^p}\frac{\tau(z_n)^2}{\tau(\varphi(z_n))^{2}}\left|1-\left|
\frac{K_{\varphi(z_n)}(\psi(z_n))}{K_{\varphi(z_n)}(\varphi(z_n))}\right|\right|^p\\
&\gtrsim\frac{\omega(z_n)^p}{\omega(\varphi(z_n))^p}\left|1-\frac{|K_{\varphi(z_n)}(\psi(z_n))|}{\|K_{\varphi(z_n)}\|\|K_{\psi(z_n)}\|}\right|^p\nonumber\\
&\gtrsim\frac{\omega(z_n)^p}{\omega(\varphi(z_n))^p}\rho_{\tau,\varphi,\psi}(z_n)^{2p}\gtrsim1\nonumber
\end{align}
as $n\rightarrow\infty$. Thus, we derive a contradiction to the compactness of $C_\varphi-C_\psi$.

\end{proof}

\begin{cor}\label{necessarybdd}
Let $\omega\in\mathcal{W}$. For $0<p<\infty$,
\begin{align*}
\|C_\varphi-C_{\psi}\|_{p}&\gtrsim\limsup_{\rho_{\tau,\varphi,\psi}(z)\rightarrow1}\left(\frac{\omega(z)}{\omega(\varphi(z))}+
\frac{\omega(z)}{\omega(\psi(z))}\right).
\end{align*}
\end{cor}
\begin{proof}
We note that $|z_n|\rightarrow1$ whenever $\lim_{n\rightarrow\infty}\rho_{\tau,\varphi,\psi}(z_n)=1$, so  $\lim_{n\rightarrow\infty}\frac{\tau(z_n)}{\tau(\varphi(z_n))}\gtrsim1$ if $\lim_{n\rightarrow\infty}\frac{\omega(z_n)}{\omega(\varphi(z_n))}\neq0$ by Lemma \ref{delight}. We also note that (\ref{goal3}) must converge to $0$ when $\lim_{n\rightarrow\infty}\frac{\omega(z_n)}{\omega(\varphi(z_n))}=0$ by Lemma \ref{delight}, thus we obtain our desired lower bound from the result of Proposition \ref{mainresult4}.

\end{proof}

By the compactness criterion for linear operators on $A^p(\omega)$, we obtain the lower bound of the essential norm of $C_\varphi-C_\psi$ promptly.
\begin{cor}
Let $\omega\in\mathcal{W}$. For $0<p<\infty$,
\begin{align*}
\|C_\varphi-C_{\psi}\|_{e}&\gtrsim\limsup_{\rho_{\tau,\varphi,\psi}(z)\rightarrow1}\left(\frac{\omega(z)}{\omega(\varphi(z))}+
\frac{\omega(z)}{\omega(\psi(z))}\right).
\end{align*}
\end{cor}
\begin{proof}
For any compact operator $K$ on $A^p(\omega)$, we have
\begin{align*}
\|C_\varphi-C_{\psi}-K\|_p&\gtrsim\limsup_{n\rightarrow\infty}|\|(C_\varphi-C_{\psi})g_n\|_p-\|Kg_n\|_p|
\\&=\limsup_{n\rightarrow\infty}\|(C_\varphi-C_{\psi})g_n\|_p,
\end{align*}
where $g_n$ is defined in (\ref{testft}) for any sequence $\{z_n\}$ with $|z_n|\rightarrow1$. The rest of the proof is the same as the proof of Corollary \ref{necessarybdd}.
\end{proof}
Now, combining Proposition \ref{mainresult3}, Proposition \ref{mainresult4} with the following Corollary, we complete the proof of Theorem \ref{k-thm1}.
\begin{lem}\cite[Lemma 3.1]{P1}\label{easypart}
Let $\omega\in\mathcal{W}$ and $0<p<\infty$. Let $u(z)$ be a nonnegative bounded measurable function on $\D$ and $C_\varphi$ be bounded on $A^p(\omega)$. If
\begin{align*}
\lim_{|z|\rightarrow1}u(z)\frac{\omega(z)}{\omega(\varphi(z))}=0
\end{align*}
then $uC_\varphi:A^p(\omega)\rightarrow L^p(\omega dA)$ is compact.
\end{lem}

\begin{cor}\label{result1}
Let $\omega\in\mathcal{W}$ and $0<p<\infty$. Suppose $C_\varphi, C_\psi$ are bounded on $A^p(\omega)$. If
\begin{align*}
\lim_{|z|\rightarrow1}\rho_{\tau,\varphi,\psi}(z)\left(\frac{\omega(z)}{\omega(\varphi(z))}+\frac{\omega(z)}{\omega(\psi(z))}\right)=0
\end{align*}
then $\rho_{\tau,\varphi,\psi} C_\varphi$ and $\rho_{\tau,\varphi,\psi} C_\psi$ are compact from $A^p(\omega)$ into $L^p(\omega dA)$.
\end{cor}

\section{Example}
In \cite{KM}, Kriete and MacCluer gave the characterization with respect to the angular derivative as follows: $\varphi$ induces an unbounded composition operator $C_\varphi$ on $A^2(\omega)$ if there is a boundary point $\zeta$ such that $|\varphi'(\zeta)|<1$. Moreover,
\begin{align}\label{angulcomp}
C_\varphi\quad \text{is compact on } A^2(\omega) \quad \iff \quad |\varphi'(\zeta)|>1,\quad \zeta\in\partial\D.
\end{align}
\indent Recall that $\varphi$ has a finite angular derivative $\varphi'(\zeta)$ at a boundary point $\zeta$ if we denote $\varphi(\zeta):=\angle\lim_{z\rightarrow\zeta}\varphi(z)$,
\begin{align*}
\varphi'(\zeta):=\angle\lim_{\substack{z\to \zeta\\ z\in \Gamma(\zeta,\alpha)}}\frac{\varphi(z)-\varphi(\zeta)}{z-\zeta}<\infty\quad\text{for each}\quad\alpha>1
\end{align*}
where $\Gamma(\zeta,\alpha)=\{z\in\D:|z-\zeta|<\alpha(1-|z|)\}$. It is well-known as the Julia-Caratheodory Theorem which tells
\begin{align}\label{Julia}
|\varphi'(\zeta)|=\liminf_{z\rightarrow\zeta}\frac{1-|\varphi(z)|}{1-|z|}<\infty.
\end{align}
Consider two analytic self-maps $\varphi$ and $\psi$ having finite
angular derivatives at $\zeta$. We say that $\varphi$ and $\psi$ have the same $M$ order data at $\zeta$ if $\varphi, \psi$ are $M$-th continuously differentiable at $\zeta$ and $\varphi^{(n)}(\zeta)=\psi^{(n)}(\zeta)$ for $n=0,1,\ldots,M$.\\
\indent For $k>0$, it is said that $\varphi(\D)$ has order of contact at most $k$ with $\partial\D$ if for each $\zeta\in\partial\D$ there is a neighborhood $\mathcal{N}(\zeta)$ centered at $\zeta$ so that
\begin{align*}
\inf\left\{\frac{1-|z|}{|\zeta-z|^k}:z\in\varphi(\D)\cap\mathcal{N}(\zeta)\right\}>0.
\end{align*}

\begin{pro}\label{omegatau}
Let $\omega=e^{-\eta}\in\mathcal{W}$ and $\varphi, \psi\in\mathcal{S}(\D)$. Suppose there exists the smallest integer $m\geq1$ satisfying $\frac{\tau(r)}{(1-r)^m}\gtrsim1$ and $\varphi, \psi$ are $M$-th continuously differentiable at $\zeta$ for $M\geq1$. If $\varphi, \psi$ have the same $M$ order data at $\zeta$ and $\varphi$ has order of contact at most $\frac{M}{m}$ at $\zeta$ then $\lim_{z\rightarrow\zeta}\rho_\tau(\varphi(z),\psi(z))=0$.
\end{pro}
\begin{proof}
Assume that $\varphi(\zeta)=\psi(\zeta)$, $\varphi^{(n)}(\zeta)=\psi^{(n)}(\zeta)$, $n=1,\ldots,M$, then by the Taylor expansions of $\varphi$ and $\psi$, we have
\begin{align}\label{Talyor}
|\varphi(z)-\psi(z)|=|h(z)|
\end{align}
where $h(z)=o(|z-\zeta|^{M})$. Therefore, we obtain
\begin{align*}
\frac{|\varphi(z)-\psi(z)|}{\tau(\varphi(z))}&\lesssim\frac{|h(z)|}{|\varphi(z)-\varphi(\zeta)|^{M}}
\frac{|\varphi(z)-\varphi(\zeta)|^{M}}{(1-|\varphi(z)|)^{m}}\\
&\lesssim\frac{|h(z)|}{|z-\zeta|^{M}}\longrightarrow0
\end{align*}
when $z\rightarrow\zeta$ since $\frac{|\varphi(z)-\varphi(\zeta)|^{M}}{(1-|\varphi(z)|)^{m}}\lesssim1$ near $\zeta$. Thus, we complete our proof.
\end{proof}
We remark that the case $m=1$ of the result above implies the standard weight case, which was proved in \cite{M}. In conjunction with Theorem \ref{k-thm1}, the following result can be obtained immediately from Proposition \ref{omegatau}.
\begin{cor}\label{deriversion}
Let $\omega=e^{-\eta}\in\mathcal{W}$ and $C_\varphi, C_\psi$ be bounded on $A^p(\omega)$. Suppose there exists the smallest integer $m\geq1$ satisfying $\frac{\tau(r)}{(1-r)^m}\gtrsim1$ and $\varphi, \psi$ are $M$-th continuously differentiable at $\zeta$ for $M\geq1$. If $\varphi, \psi$ have the same $M$ order data at $\zeta$ and $\varphi$ has order of contact at most $\frac{M}{m}$ at $\zeta$ then $C_\varphi-C_\psi$ is compact on $A^p(\omega)$.
\end{cor}

Finally, we provide an example showing that $C_\varphi$ and $C_\psi$ are not compact on $A^2(\omega)$ but their difference is compact.
\begin{ex}
Consider the Bergman space having the weight $\omega(z)=e^{-\frac{1}{1-|z|}}$. Put
\begin{align*}
\varphi(z)=\frac{1+z^2}{2}\quad\text{and}\quad\psi(z)=\varphi(z)+\epsilon(1-z^2)^5,\quad0<\epsilon<\frac{1}{2^7}.
\end{align*}
Then $C_\varphi$ and $C_\psi$ are not compact on $A^2(\omega)$ but $C_\varphi-C_\psi$ is compact.
\end{ex}
\begin{proof}
First, we note that the following estimates hold:
\begin{align}\label{calculation}
\psi(z)=\varphi(z)+32\epsilon(1-\varphi(z))^5\quad\text{and}\quad1-|\varphi(z)|^2\geq|1-\varphi(z)|^2,\quad z\in\D.
\end{align}
Thus, $\psi\in S(\D)$ since
\begin{align}\label{lesson}
|\psi(z)|\leq|\varphi(z)|+32\epsilon(1-|\varphi(z)|^2)^{2}\leq1.
\end{align}
Here, note that the second equality of (\ref{lesson}) holds only when $z=1,-1$. In order to see the boundedness of $C_\varphi$ and $C_\psi$, we calculate
\begin{align*}
\lim_{z\rightarrow\zeta\neq1,-1}\frac{\omega(z)}{\omega(\varphi(z))}=\lim_{z\rightarrow\zeta\neq1,-1}\exp\left(-\frac{1}{1-|z|}+
\frac{1}{1-\big|\frac{1+z^2}{2}\big|}\right)=0.
\end{align*}
In the same way, $\lim_{z\rightarrow\zeta\neq1,-1}\frac{\omega(z)}{\omega(\psi(z))}=0$ by (\ref{lesson}). Moreover, since
\begin{align}\label{calculation2}
\frac{|\varphi(z)|-|z|}{(1-|z|)(1-|\varphi(z)|)}\leq\frac{1}{2}\frac{(1-|z|)^2}{(1-|z|)(1-|\varphi(z)|)}\leq\frac{1}{2}\frac{1-|z|}{1-|\varphi(z)|},
\end{align}
and $|\varphi'(1)|=|\varphi'(-1)|=1=\limsup_{z\rightarrow1,-1}\frac{1-|z|}{1-|\varphi(z)|}$ by (\ref{Julia}), we have
\begin{align*}
\limsup_{z\rightarrow1,-1}\frac{\omega(z)}{\omega(\varphi(z))}&=\limsup_{z\rightarrow1,-1}
\exp\left(\frac{|\varphi(z)|-|z|}{(1-|z|)(1-|\varphi(z)|)}\right)\\&\leq
\limsup_{z\rightarrow1,-1}\exp\left(\frac{1}{2}\frac{1-|z|}{1-|\varphi(z)|}\right)<\infty.
\end{align*}
Likewise, using (\ref{lesson}) and (\ref{calculation2}), we have
\begin{align*}
\frac{|\psi(z)|-|z|}{(1-|z|)(1-|\psi(z)|)}&\leq\frac{|\varphi(z)|-|z|+4(1-|z|+|z|-|\varphi(z)|)^2}{(1-|z|)(1-|\psi(z)|)}\\
&\leq\frac{1/2(1-|z|)^2+8(1-|z|)^2+8(|\varphi(z)|-|z|)^2}{(1-|z|)(1-|\psi(z)|)}\\
&\leq\frac{11(1-|z|)^2}{(1-|z|)(1-|\psi(z)|)}\lesssim\frac{1-|z|}{1-|\psi(z)|}.
\end{align*}
Thus, $\limsup_{z\rightarrow1,-1}\frac{\omega(z)}{\omega(\psi(z))}<\infty$ since $|\psi'(1)|=|\psi'(-1)|=1$. Therefore, we confirm that $\varphi$ and $\psi$ induce the bounded composition operators but they are not compact on $A^2(\omega)$ from (\ref{bddwithoutU}) and (\ref{angulcomp}). Now, it is easily checked that $\varphi, \psi$ have the same $4$-order data at $1,-1$ and $\varphi$ has order of contact at most $2$ from (\ref{calculation}) since
\begin{align*}
\frac{1-\big|\varphi(z)\big|}{\big|1-\varphi(z)|^2}\geq\frac{1}{2}\quad\text{near}\quad z=1, -1.
\end{align*}
Letting $\tau(z)=(1-|z|)^{3/2}$, we obtain $\rho_{\tau,\varphi,\psi}(z)\rightarrow0$ when $z\rightarrow1,-1$ if we take $m=2$ and $M=4$ in Proposition \ref{omegatau}. Thus, we conclude $C_\varphi-C_\psi$ is compact by Theorem \ref{k-thm1}.

\end{proof}

\section{Hilbert-Schmidt Difference}
Let $H_1$ and $H_2$ be Hilbert spaces. Recall that a bounded linear operator $T:H_1\rightarrow H_2$ is Hilbert-Schmidt if
\begin{align*}
\|T\|_{HS}=\sum_{n=0}^{\infty}\|Te_n\|^2<\infty
\end{align*}
where $\{e_n\}$ is an arbitrary orthonormal basis for the Hilbert space $H_1$. Thus, using the definition of the reproducing kernel (\ref{definition}), we easily obtain the Hilbert-Schmidt norm of $uC_\varphi$ on $A^2(\omega)$,
\begin{align}\label{defhs}
\|uC_\varphi\|_{HS}&=\sum_{n=1}^{\infty}\int_{\D}|u(z)|^2|e_n(\varphi(z))|^2\omega(z)^2dA\nonumber\\
&=\int_\D|u(z)|^2\|K_{\varphi(z)}\|^2\omega(z)^2dA(z).
\end{align}

\begin{lem}\label{normHS}
Let $\omega\in\mathcal{W}$ and $u$ be a measurable function on $\D$. Then
\begin{align*}
\left\|u(C_{\varphi}-C_\psi)\right\|^2_{HS}=\int_\D |u(z)|^2\|K_{\varphi(z)}-K_{\psi(z)}\|^2\omega(z)^2dA(z).
\end{align*}
\end{lem}
\begin{proof}
Consider an arbitrary orthonormal basis $\{e_n\}$ for $A^2(\omega)$. By the definition of the Hilbert-Schmidt norm and the reproducing kernel of $A^2(\omega)$, we obtain that
\begin{align*}
&\sum_{n=0}^{\infty}\|(uC_\varphi-uC_\psi)e_n\|^2\\&=\sum_{n=0}^{\infty}[\|ue_n(\varphi)\|^2+\|ue_n(\psi)\|^2-2\re\langle ue_n(\varphi),ue_n(\psi)\rangle_\omega]\\
&=\int_\D|u(z)|^2\sum_{n=0}^{\infty}\big[|e_n(\varphi(z))|^2+|e_n(\psi(z))|^2-2\re e_n(\varphi(z))\overline{e_n(\psi(z))}\big]\omega(z)^2dA(z)\\
&=\int_\D|u(z)|^2(\|K_{\varphi(z)}\|^2+\|K_{\psi(z)}\|^2-2\re K(\varphi(z),\psi(z)))\omega(z)^2dA(z).
\end{align*}
\end{proof}

Using Lemma \ref{normHS} and Theorem \ref{diffkernelnorm}, we obtain a Hilbert-Schmidt norm estimate for $C_\varphi-C_\psi$ on $A^2(\omega)$ involving the distance $\rho_{\tau,\varphi,\psi}$:
\begin{pro}\label{key3}
Let $\omega=e^{-\eta}\in\mathcal{W}$. Then
\begin{align*}
\left\|C_{\varphi}-C_\psi\right\|^2_{HS}\approx\int_\D\rho_{\tau,\varphi,\psi}(z)^2(\|K_{\varphi(z)}\|^2+\|K_{\psi(z)}\|^2)
\omega(z)^2dA(z).
\end{align*}
\end{pro}

Now, from (\ref{defhs}) and Proposition \ref{key3}, we obtain the following result promptly.
\begin{thm}\label{mainresult5}
Let $\omega\in\mathcal{W}$. Then $\rho_{\tau,\varphi,\psi} C_\varphi$, $\rho_{\tau,\varphi,\psi} C_\psi$ are Hilbert-Schmidt operators from $A^2(\omega)$ into $L^2(\omega dA)$ if and only if $C_\varphi-C_\psi$ is a Hilbert-Schmidt operator on $A^2(\omega)$.
\end{thm}

Finally, we close this section with observing the path components of $\mathcal{C}_{HS}(A^2(\omega))$ which denotes the space of all composition operators on $A^2(\omega)$ endowed with topology induced by the metric:
\begin{align*}
d(C_\varphi,C_\psi):=\left\{
 \begin{array}{cc}
\frac{\|C_\varphi-C_\psi\|_{HS}}{1+\|C_\varphi-C_\psi\|_{HS}}, &  \|C_\varphi-C_\psi\|_{HS}<\infty\\
1,&  \|C_\varphi-C_\psi\|_{HS}=\infty.\\
 \end{array}
 \right.
\end{align*}
It is said that $C_\varphi$, $C_\psi$ are in the same path component of $\mathcal{C}_{HS}(A^2(\omega))$ if there exists a continuous path $C_{\gamma(s)}$ with respect to $s$ in $\mathcal{C}_{HS}(A^2(\omega))$ such that $\gamma(0)=\varphi$, $\gamma(1)=\psi$. Define the set
\begin{align*}
U(C_\varphi):=\{C_\psi:\|C_\varphi-C_\psi\|_{HS}<\infty\}.
\end{align*}
In the same sense, $U(C_\varphi)$ is said to be convex if $C_{(1-s)\varphi+s\psi}$ is continuous with respect to $s$ in $\mathcal{C}_{HS}(A^2(\omega))$ for every $C_\psi\in U(C_\varphi)$. Now, we will prove that the set $U(C_\varphi)$ is the path connected component of $C_\varphi$ in $\mathcal{C}_{HS}(A^2(\omega))$ and is also convex. \cite{HJM, CHK} contain some results for the connected component of the space of composition operators acting on the Hardy and the standard weighted Bergman spaces under the Hilbert-Schmidt norm topologies, respectively.
\begin{lem}\label{careful}
Let $\omega\in\mathcal{W}$. For $z,w\in\D$, denote $z_s:=(1-s)z+sw$ for $s\in[0,1]$. Then
\begin{align*}
\rho_\tau(z_s,z_t)\leq C\rho_\tau(z,w)
\end{align*}
where $C>0$ is independent of $s,t$.
\end{lem}
\begin{proof}
We first consider the case $d_\tau(z,w)<R$. By Lemma \ref{setinclus}, there exists $0<C_1<1$ satisfying
\begin{align}\label{unifbdd}
d_\tau(z_s,z_t)&\leq\frac{|z_t-z_s|}{\min(\tau(z_s),\tau(z_t))}\leq\frac{|s-t||z-w|}{\min(\tau(z),\tau(w))}\leq\frac{1}{C_1}d_\tau(z,w).
\end{align}
Moreover, it is clear that $\rho_\tau(z_s,z_t)\leq \frac{1}{1-e^{-R}}\rho_\tau(z,w)$ for $d_\tau(z,w)\geq R$.
\end{proof}
\begin{thm}\label{mainresult1}
Let $\omega\in\mathcal{W}$ and $C_\varphi, C_\psi$ be bounded on $A^2(\omega)$. Then $C_\varphi-C_\psi$ is Hilbert-Schmidt on $A^2(\omega)$ if and only if $C_\varphi$ and $C_\psi$ lie in the same path component of $\mathcal{C}_{HS}(A^2(\omega))$.
\end{thm}
\begin{proof}
Assume that there exists a continuous path $C_{\gamma(s)}:[0,1]\rightarrow \mathcal{C}_{HS}(A^2(\omega))$. Then $C_{\gamma(s)}$ is uniformly continuous on $[0,1]$ so that given $\epsilon>0$, there is a partition $\{s_0=0,s_1,\ldots,s_{N-1},s_N=1\}\subset[0,1]$ such that $\|C_{\gamma(s_i)}-C_{\gamma(s_{i+1})}\|_{HS}<\epsilon$ for all $i=0,1,\ldots,N-1$. Thus $C_\varphi-C_\psi$ is Hilbert-Schmidt on $A^2(\omega)$ by the triangle inequality. To show the necessity, we will show that
\begin{align*}
\lim_{t\rightarrow s}\|C_{\varphi_s}-C_{\varphi_t}\|_{HS}=0
\end{align*}
when $\varphi_s=(1-s)\varphi+s\psi$ for $0\leq s\leq1$. Since $|\varphi_s(z)|\leq\max(|\varphi(z)|,|\psi(z)|)$ and $\|K_z\|$ increases with $|z|$ by (\ref{kernelesti}), we have $\|K_{\varphi_s(z)}\|\lesssim\|K_{\varphi(z)}\|+\|K_{\psi(z)}\|$. Therefore, using the results of Lemma \ref{normHS}, Proposition \ref{key3} and Lemma \ref{careful}, we have
\begin{align}\label{goal4}
\|C_{\varphi_s}-C_{\varphi_t}\|_{HS}&\lesssim\int_{\D}\rho_{\tau,\varphi_s,\varphi_t}(z)^2(\|K_{\varphi_s(z)}\|^2+\|K_{\varphi_t(z)}\|^2)\omega(z)^2dA(z)\\
&\lesssim\int_{\D}\rho_{\tau,\varphi,\psi}(z)^2(\|K_{\varphi(z)}\|^2+\|K_{\psi(z)}\|^2)\omega(z)^2dA(z)<\infty\nonumber.
\end{align}
By the Dominated Convergence Theorem, (\ref{goal4}) vanishes when $t\rightarrow s$ for $\rho_{\tau,\varphi_s,\varphi_t}(z)\rightarrow0$.
\end{proof}
The inequality above gives the following result immediately.
\begin{cor}
Let $\omega\in\mathcal{W}$ and $C_\varphi, C_\psi$ be bounded on $A^2(\omega)$. Define $\varphi_s=(1-s)\varphi+s\psi$ for $0\leq s\leq1$. Then $C_\varphi-C_\psi$ is Hilbert-Schmidt on $A^2(\omega)$ if and only if $C_{\varphi_s}-C_{\varphi_t}$ is Hilbert-Schmidt on $A^2(\omega)$ for $0\leq s,t\leq1$.
\end{cor}

\bibliographystyle{amsplain}

\end{document}